\documentclass[12pt]{amsart}
\usepackage{amscd,amsmath,amssymb,amsfonts}
\theoremstyle{plain}
\newtheorem{thm}{Theorem}
\newtheorem{lem}[thm]{Lemma}

\theoremstyle{definition}
\newtheorem{defn}[thm]{Definition}

\numberwithin{thm}{section} \numberwithin{equation}{section}

\newcommand{\ga}[2]{\begin{gather}\label{#1}#2 \end{gather}}

\newcommand{\surj}{\twoheadrightarrow}
\newcommand{\inj}{\hookrightarrow}

\newcommand{\Spec}{{\rm Spec \,}}

\newcommand{\sE}{{\mathcal E}}

\newcommand{\sL}{{\mathcal L}}

\newcommand{\sO}{{\mathcal O}}

\newcommand{\sV}{{\mathcal V}}

\begin{document}

\title{Stability of direct images under Frobenius morphism}
\author{Xiaotao Sun}
\address{Chinese Academy of Mathematics and Systems Science, Beijing, P. R. of China}
\email{xsun@math.ac.cn}
\address{}
\date{August 18, 2006}
\thanks{Partially supported by the DFG Leibniz Preis of Esnault-Viehweg}
\begin{abstract} Let $X$ be a smooth projective variety over an
algebraically field $k$ with ${\rm char}(k)=p>0$ and $F:X\to X_1$
be the relative Frobenius morphism. When ${\rm dim}(X)=1$, we
prove that $F_*W$ is a stable bundle for any stable bundle $W$
(Theorem \ref{thm1.3}). As a step to study the question for higher
dimensional $X$, we generalize the canonical filtration (defined
by Joshi-Ramanan-Xia-Yu for curves) to higher dimensional $X$
(Theorem \ref{thm2.6}).
\end{abstract}
\maketitle
\begin{quote}

\end{quote}
\section{Introduction}

Let $X$ be a smooth projective variety over an algebraically field
$k$ with ${\rm char}(k)=p>0$ and $F:X\to X_1$ be the relative
Frobenius morphism. When ${\rm dim}(X)=1$, Lange and Pauly proved
that $F_*\sL$ is a stable bundle for a line bundle $\sL$ (cf.
\cite[Proposition 1.2 ]{LP}). The first result in this paper is
that stability of $W$ implies stability of $F_*W$.

Recall that for a Galois {\'e}tale $G$-cover $f:Y\to X$ and a
semistable bundle $W$ on $Y$, to prove semistability of $f_*W$,
one uses the fact that $f^*(f_*W)$ decomposes into pieces of
$W^{\sigma}$ ($\sigma\in G$). To imitate this idea for $F:X\to
X_1$, we need a similar decomposition of $V=F^*(F_*W)$. Indeed,
use the canonical connection $\nabla: V\to V\otimes\Omega^1_X$,
Joshi-Ramanan-Xia-Yu have defined in \cite{JRXY} for ${\rm
dim}(X)=1$ a canonical filtration
$$0=V_p\subset V_{p-1}\subset\cdots\subset V_i\subset
V_{i-1}\subset\cdots V_1\subset V_0=V$$ such that
$V_i/V_{i+1}\cong W\otimes(\Omega^1_X)^{\otimes i}$. For any
$0\neq\sE\subset F_*W$, let $$0\subset V_m\cap
F^*\sE\subset\,\cdots\,\subset V_1\cap F^*\sE\subset V_0\cap
F^*\sE=F^*\sE$$ be the induced filtration. Then we can show (cf.
Lemma \ref{lem1.2})
$$\frac{V_{i-1}\cap F^*\sE}{V_i\cap F^*\sE}\neq 0 \quad {\rm for}
\quad 1\le i\le m+1.$$ Using the induced filtration and stability
of $W\otimes (\Omega^1_X)^{\otimes i}$, we have
$$\mu(F_*W)-\mu(\sE)\ge \frac{g-1}{p}\left(p-1-\frac{2}{{\rm
rk}(\sE)}\cdot\sum_{i=1}^{m+1}(i-1){\rm rk}(\frac{V_{i-1}\cap
F^*\sE}{V_i\cap F^*\sE})\right).$$ When $W$ is a line bundle, all
$\frac{V_{i-1}\cap F^*\sE}{V_i\cap F^*\sE}$ must be line bundles
and ${\rm rk}(\sE)=m+1$. Then above inequality implies the
stability of $F_*W$ immediately. For higher rank bundles $W$, we
need more analysis of the rank of $\frac{V_{i-1}\cap
F^*\sE}{V_i\cap F^*\sE}$.

It is a natural question to study $F_*W$ for ${\rm dim}(X)=n>1$.
As the first step, we generalize the canonical filtration to
higher dimensional $X$. Its definition can be generalized
straightforwardly by using the canonical connection $\nabla:V\to
V\otimes\Omega^1_X$. The second result of this paper is that there
exists a canonical filtration
$$0=V_{n(p-1)+1}\subset V_{n(p-1)}\subset\cdots\subset V_1\subset
V_0=V=F^*(F_*W)$$ such that $\nabla$ induces $V_i/V_{i+1}\cong
W\otimes (\Omega^1_X)^{[i]}$, where $(\Omega^1_X)^{[i]}\subset
(\Omega^1_X)^{\otimes i}$ is a subbundle given by a representation
of ${\rm GL}(n)$ (cf. Definition \ref{defn2.4}). In characteristic
zero, $(\Omega^1_X)^{[i]}={\rm Sym}^i(\Omega^1_X)$. In
characteristic $p>0$, we have $(\Omega^1_X)^{[i]}\cong {\rm
Sym}^i(\Omega^1_X)$ for $i<p$. The general question would be: how
to bound the instability of $F_*W$ by instability of
$W\otimes(\Omega^1_X)^{[i]}$ ?

When I was preparing the last section of this paper, Mehta and
Pauly posted a preprint \cite{MP}, in which they prove, in a
different mothed, that semistability of $W$ implies semistability
of $F_*W$. But they do not prove that stability of $W$ implies
stability of $F_*W$.

\section{The case of curves}
Let $k$ be an algebraically closed field of characteristic $p>0$
and $X$ be a smooth projective curve over $k$. Let $F:X\to X_1$ be
the relative $k$-linear Frobenius morphism, where
$X_1:=X\times_kk$ is the base change of $X/k$ under the Frobenius
$\Spec(k)\to\Spec(k)$. Let $W$ be a vector bundle on $X$ and
$V=F^*(F_*W)$. It is known (\cite[Theorem 5.1]{K}) that $V$ has an
$p$-curvature zero connection $\nabla: V\to V\otimes\Omega^1_X$.
In \cite[Section 5]{JRXY}, the authors defined a canonical
filtration \ga{1.1}{0=V_p\subset V_{p-1}\subset\cdots\subset
V_i\subset V_{i-1}\subset\cdots V_1\subset V_0=V} where $V_1={\rm
ker}(V=F^*F_*W\surj W)$ and \ga{1.2}{V_{i+1}={\rm
ker}(V_i\xrightarrow{\nabla} V\otimes\Omega^1_X\to
V/V_i\otimes\Omega^1_X).} The following lemma belongs to them (cf.
\cite[Theorem 5.3]{JRXY}).
\begin{lem}\label{lem1.1}
\begin{itemize}\item[(i)] $V_0/V_1\cong W$, $\nabla(V_{i+1})\subset V_i\otimes\Omega^1_X$
for $i\ge 1$.
\item[(ii)]
$V_i/V_{i+1}\xrightarrow{\nabla}(V_{i-1}/V_i)\otimes\Omega^1_X$ is
an isomorphism for $1\le i\le p-1$.
\item[(iii)] If $g\ge 2$ and $W$ is semistable, then the canonical
filtration \eqref{1.1} is nothing but the Harder-Narasimhan
filtration.
\end{itemize}
\end{lem}
\begin{proof} (i) follows by the definition, which and (ii) imply
(iii). To prove (ii), let $I_0=F^*F_*\sO_X$, $I_1={\rm
ker}(F^*F_*\sO_X\surj \sO_X)$ and
 \ga{1.3} {I_{i+1}={\rm ker}(I_i\xrightarrow{\nabla}
I_0\otimes\Omega^1_X\surj I_0/I_i\otimes\Omega^1_X)} which is the
canonical filtration \eqref{1.1} in the case $W=\sO_X$.

(ii) is clearly a local problem, we can assume $X=\Spec(k[[x]])$
and $W=k[[x]]^{\oplus r}$. Then $V_0:=V=F^*(F_*W)=I_0^{\oplus r}$,
$V_i=I_i^{\oplus r}$ and
\ga{1.4}{V_i/V_{i+1}=(I_i/I_{i+1})^{\oplus
r}\xrightarrow{\oplus\nabla}
(I_{i-1}/I_i\otimes\Omega^1_X)^{\oplus
r}=V_{i-1}/V_i\otimes\Omega^1_X.} Thus it is enough to show that
\ga{1.5}{I_i/I_{i+1}\xrightarrow{\nabla}
(I_{i-1}/I_i)\otimes\Omega^1_X} is an isomorphism. Locally,
$I_0=k[[x]]\otimes_{k[[x^p]]}k[[x]]$ and
\ga{1.6}{\nabla:k[[x]]\otimes_{k[[x^p]]}k[[x]]\to
I_0\otimes_{\sO_X}\Omega^1_X,} where $\nabla(g\otimes f)=g\otimes
f'\otimes{\rm d}x$. The $\sO_X$-module \ga{1.7} {I_1:={\rm
ker}(k[[x]]\otimes_{k[[x^p]]}k[[x]]\surj k[[x]])} has a basis $\{
x^k\otimes 1-1\otimes x^k\}_{1\le k\le p-1}$. Notice that $I_1$ is
also an ideal of the $\sO_X$-algebra
$I_0=k[[x]]\otimes_{k[[x^p]]}k[[x]]$, let $\alpha=x\otimes
1-1\otimes x$, then $\alpha^k\in I_1$. It is easy to see that
$\alpha,\,\alpha^2,\,\ldots,\,\alpha^{p-1}$ is a basis of the
$\sO_X$-module $I_1$ (notice that $\alpha^p=x^p\otimes 1-1\otimes
x^p=0$), and \ga{1.8} {\nabla(\alpha^k)=-k\alpha^{k-1}\otimes {\rm
d}x.} Thus, as a free $\sO_X$-module, $I_i$ has a basis
$\{\alpha^i,\,\alpha^{i+1},\ldots,\,\alpha^{p-1}\}$, which means
that $I_i/I_{i+1}$ has a basis $\alpha^i$,
$(I_{i-1}/I_i)\otimes\Omega^1_X$ has a basis
$\alpha^{i-1}\otimes{\rm d}x$ and
$\nabla(\alpha^i)=-i\alpha^{i-1}\otimes{\rm d}x$. Therefore
$\nabla$ induces the isomorphism \eqref{1.5} since $(i,p)=1$,
which implies the isomorphism in (ii).
\end{proof}

\begin{lem}\label{lem1.2} Let $\sE\subset F_*W$ be a nontrivial
subsheaf and let \ga{1.9} {0\subset V_m\cap
F^*\sE\subset\,\cdots\,\subset V_1\cap F^*\sE\subset V_0\cap
F^*\sE=F^*\sE} be the induced filtration. Then \ga{1.10}
{\frac{V_{i-1}\cap F^*\sE}{V_i\cap F^*\sE}\neq 0 \quad {\rm for}
\quad 1\le i\le m+1.}
\end{lem}
\begin{proof} Firstly, by adjunction formula, $F^*\sE\inj
V=F^*(F_*W)\surj W$ is nontrivial. Thus $V_0\cap F^*\sE/V_1\cap
F^*\sE$ is nontrivial. On the other hand, for any $i\ge 2$, the
morphism $V_{i-1}\cap F^*\sE\inj V=F^*(F_*W)\surj W$ is trivial,
which implies, by adjunction formula, that there is no subsheaf
$j:\sE'\inj F_*W$ such that $F^*j:V_{i-1}\cap F^*\sE\cong
F^*\sE'\inj V$ is the inclusion $V_{i-1}\cap F^*\sE\inj V$.
However, by the definition of canonical filtration \eqref{1.1},
$V_{i-1}\cap F^*\sE=V_i\cap F^*\sE$ implies that \ga{1.11}
{\nabla(V_{i-1}\cap F^*\sE)\subset (V_{i-1}\cap
F^*\sE)\otimes\Omega^1_X.} By \cite[Theorem 5.1]{K}, this means
that there is an $j:\sE'\inj F_*W$ such that $F^*j:V_{i-1}\cap
F^*\sE\cong F^*\sE'\inj V$ is the inclusion $V_{i-1}\cap
F^*\sE\inj V$. We get contradiction.
\end{proof}

\begin{thm}\label{thm1.3} If $W$ is a stable vector bundle, then
$F_*W$ is a stable vector bundle. In particular, if $W$ is
semistable, then $F_*W$ is semistable.
\end{thm}
\begin{proof} Let $\sE\subset F_*W$ be a nontrivial subbundle and
\ga{1.12}{0\subset V_m\cap F^*\sE\subset\,\cdots\,\subset V_1\cap
F^*\sE\subset V_0\cap F^*\sE=F^*\sE} be the induced filtration.
Let $r_{i-1}={\rm rk}(\frac{V_{i-1}\cap F^*\sE}{V_i\cap F^*\sE})$
be the ranks of quotients. Then, by the filtration \eqref{1.12},
we have \ga{1.13}{\mu(F^*\sE)=\frac{1}{{\rm
rk}(F^*\sE)}\sum^{m+1}_{i=1}r_{i-1}\mu(\frac{V_{i-1}\cap
F^*\sE}{V_i\cap F^*\sE}).} By Lemma \ref{lem1.1},
$V_{i-1}/V_i\cong W\otimes(\Omega^1_X)^{\otimes(i-1)}$ is stable,
we have \ga{1.14} {\mu(\frac{V_{i-1}\cap F^*\sE}{V_i\cap
F^*\sE})\le\mu(W)+2(g-1)(i-1).} Then, notice that
$\mu(V)=\mu(W)+(p-1)(g-1)$, we have
\ga{1.15}{\mu(F_*W)-\mu(\sE)\ge \frac{2g-2}{p\cdot{\rm
rk}(\sE)}\cdot\sum^{m+1}_{i=1}(\frac{p+1}{2}-i)r_{i-1}} which
becomes into an equality if and only if the inequalities in
\eqref{1.14} become into equalities.

It is clear by \eqref{1.15} that $\mu(F_*W)-\mu(\sE)>0$ if $m\le
\frac{p-1}{2}$. Thus we assume that $m>\frac{p-1}{2}$. On the
other hand, since the isomorphisms
$V_i/V_{i+1}\xrightarrow{\nabla}(V_{i-1}/V_i)\otimes\Omega^1_X$ in
Lemma \ref{lem1.1} (ii) induce the injections
 \ga{1.16} {\frac{V_i\cap F^*\sE}{V_{i+1}\cap F^*\sE}\inj
 \frac{V_{i-1}\cap F^*\sE}{V_i\cap F^*\sE}\otimes\Omega^1_X}
we have \ga{1.17}{r_0\ge r_1\ge \cdots\ge r_{i-1}\ge
r_i\ge\cdots\ge r_m.} Then, when $m>\frac{p-1}{2}$, we can write
\ga{1.18} {\sum^{m+1}_{i=1}(\frac{p+1}{2}-i)r_{i-1}
=\sum^{\frac{p-1}{2}}_{i=1}i\cdot
r_{\frac{p-1}{2}-i}-\sum^{m-\frac{p-1}{2}}_{i=1}i\cdot
r_{\frac{p-1}{2}+i}} Note that $m\le p-1$, use \eqref{1.17} and
\eqref{1.18}, we have
\ga{1.19}{\sum^{m+1}_{i=1}(\frac{p+1}{2}-i)r_{i-1} \ge
\sum^{m-\frac{p-1}{2}}_{i=1}i\cdot(r_{\frac{p-1}{2}-i}-
r_{\frac{p-1}{2}+i})\ge 0.} Thus we always have
 \ga{1.20}{\mu(F_*W)-\mu(\sE)\ge
\frac{2g-2}{p\cdot{\rm
rk}(\sE)}\cdot\sum^{m+1}_{i=1}(\frac{p+1}{2}-i)r_{i-1}\ge 0.} If
$\mu(F_*W)-\mu(\sE)=0$, then \eqref{1.15} and \eqref{1.19} become
into equalities. That \eqref{1.15} becomes into an equality
implies inequalities in \eqref{1.14} become into equalities, which
means $r_0=r_1=\cdots=r_m={\rm rk}(W)$. Then that \eqref{1.19}
become into equalities implies $m=p-1$. Altogether imply
$\sE=F_*W$, we get contradiction. Hence $F_*W$ is a stable vector
bundle whenever $W$ is stable.
\end{proof}

\section{Generalizations to higher dimension varieties}

Let $X$ be a smooth projective variety over $k$ of dimension $n$
and $F:X\to X_1$ be the relative $k$-linear Frobenius morphism,
where $X_1:=X\times_kk$ is the base change of $X/k$ under the
Frobenius $\Spec(k)\to\Spec(k)$. Let $W$ be a vector bundle on $X$
and $V=F^*(F_*W)$. We have the straightforward generalization of
the canonical filtration to higher dimensional varieties.

\begin{defn}\label{defn2.1} Let $V_0:=V=F^*(F_*W)$,
$V_1=\ker(F^*(F_*W)\surj W)$ \ga{2.1}
{V_{i+1}:=\ker(V_i\xrightarrow{\nabla} V\otimes_{\sO_X}
\Omega^1_X\to (V/V_i)\otimes_{\sO_X}\Omega^1_X)} where $\nabla:
V\to V\otimes_{\sO_X} \Omega^1_X$ is the canonical connection (cf.
\cite[Theorem]{K}).
\end{defn}

We first consider the special case $W=\sO_X$ and give some local
descriptions. Let $I_0=F^*(F_*\sO_X)$, $I_1=\ker(F^*F_*\sO_X\surj
\sO_X)$ and \ga{2.2}
{I_{i+1}=\ker(I_i\xrightarrow{\nabla}I_0\otimes_{\sO_X}\Omega^1_X\to
I_0/I_i\otimes_{\sO_X}\Omega^1_X).}

Locally, let $X=\Spec(A)$, $I_0=A\otimes_{A^p} A$, where $
A=k[[x_1,\cdots,x_n]]$, $A^p=k[[x^p_1,\cdots,x^p_n]]$. Then the
canonical connection $\nabla: I_0\to I_0\otimes\Omega^1_X$ is
locally defined by \ga{2.3} {\nabla(g\otimes_{A^p} f)=\sum_{i=1}^n
(g\otimes_{A^p}\frac{\partial f}{\partial x_i})\otimes_A {\rm
d}x_i}  Notice that $I_0$ has an $A$-algebra structure such that
$I_0=A\otimes_{A^p}A\surj A$ is a homomorphism of $A$-algebras,
its kernel $I_1$ contains elements \ga{2.4}
{\alpha_1^{k_1}\alpha_2^{k_2}\cdots \alpha_n^{k_n},\,\,\,
\text{where $\alpha_i=x_i\otimes_{A^p} 1-1\otimes_{A^p} x_i$,\,\,
$\sum^n_{i=1} k_i\ge 1.$}} Since
$\alpha_i^p=x_i^p\otimes_{A^p}1-1\otimes_{A^p}x_i^p=0$, the set
$\{\alpha_1^{k_1}\cdots \alpha_n^{k_n}\,|\, k_1+\cdots+k_n\ge 1\}$
has $p^n-1$ elements. In fact, we have
\begin{lem}\label{lem2.2} Locally, as
free $A$-modules, we have, for all $i\ge 1$, \ga{2.5}
{I_i=\bigoplus_{k_1+\cdots+k_n\ge
i}(\alpha_1^{k_1}\cdots\alpha_n^{k_n})A.}
\end{lem}
\begin{proof} We first prove for $i=1$ that
$\{\alpha_1^{k_1}\cdots\alpha_n^{k_n}\,|\,k_1+\cdots+k_n\ge 1\}$
is a basis of $I_1$ locally. By definition, $I_1$ is locally free
of rank $p^n-1$, thus it is enough to show that as an $A$-module
$I_1$ is generated locally by $\{\alpha_1^{k_1}\cdots
\alpha_n^{k_n}\,|\, k_1+\cdots+k_n\ge 1\}$ since it has exactly
$p^n-1$ elements.

It is easy to see that as an $A$-module $I_1$ is locally generated
by $\{x_1^{k_1}\cdots
x_n^{k_n}\otimes_{A^p}1-1\otimes_{A^p}x_1^{k_1}\cdots
x_n^{k_n}\,|\,k_1+\cdots+k_n\ge 1,\,\,0\le k_i\le p-1\,\}$. It is
enough to show any $x_1^{k_1}\cdots
x_n^{k_n}\otimes_{A^p}1-1\otimes_{A^p}x_1^{k_1}\cdots x_n^{k_n}$
is a linear combination of
$\{\alpha_1^{k_1}\cdots\alpha_n^{k_n}\,|\,k_1+\cdots+k_n\ge 1\}.$
The claim is obvious when $k_1+\cdots+k_n=1$, we consider the case
$k_1+\cdots+k_n>1$. Without loss generality, assume $k_n\ge 1$ and
there are $f_{j_1,\ldots,j_n}\in A$ such that
$$x_1^{k_1}\cdots
x_n^{k_n-1}\otimes_{A^p}1-1\otimes_{A^p}x_1^{k_1}\cdots
x_n^{k_n-1}=\sum_{j_1+\cdots+j_n\ge
1}(\alpha_1^{j_1}\cdots\alpha_n^{j_n})\cdot f_{j_1,\ldots,j_n}.$$
Then we have
$$\aligned &x_1^{k_1}\cdots
x_n^{k_n}\otimes_{A^p}1-1\otimes_{A^p}x_1^{k_1}\cdots
x_n^{k_n}=\sum_{j_1+\cdots+j_n\ge
1}(\alpha_1^{j_1}\cdots\alpha_n^{j_n+1})\cdot
f_{j_1,\ldots,j_n}\\& +\sum_{j_1+\cdots+j_n\ge
1}(\alpha_1^{j_1}\cdots\alpha_n^{j_n})\cdot f_{j_1,\ldots,j_n}x_n
\quad+\quad\alpha_n\cdot (x_1^{k_1}\cdots x_n^{k_n-1}).
\endaligned$$

For $i>1$, to prove the lemma, we first show
\ga{2.6}{\nabla(\alpha_1^{k_1}\cdots\alpha_n^{k_n})=-\sum^n_{i=1}k_i
(\alpha_1^{k_1}\cdots\alpha_i^{k_i-1}\cdots\alpha_n^{k_n})
\otimes_A{\rm d}x_i}  Indeed, \eqref{2.6} is true when
$k_1+\cdots+k_n=1$. If $k_1+\cdots+k_n>1$, we assume $k_n\ge 1$
and $\alpha_1^{k_1}\cdots\alpha_n^{k_n-1}=\sum
g_j\otimes_{A^p}f_j$. Then
$$\alpha_1^{k_1}\cdots\alpha_n^{k_n}=\sum_j
x_ng_j\otimes_{A^p}f_j-\sum_jg_j\otimes_{A^p}f_jx_n\,.$$ Use
\eqref{2.3}, straightforward computations show
$$\nabla(\alpha_1^{k_1}\cdots\alpha_n^{k_n})=\alpha_n\nabla(\alpha_1^{k_1}\cdots\alpha_n^{k_n-1})-
(\alpha_1^{k_1}\cdots\alpha_n^{k_n-1})\otimes_A{\rm d}x_n$$ which
implies \eqref{2.6}. Now we can assume the lemma is true for
$I_{i-1}$ and recall that
$I_i=\ker(I_{i-1}\xrightarrow{\nabla}I_0\otimes_A\Omega^1_X\surj
(I_0/I_{i-1})\otimes_A\Omega^1_X)$. For any
$$\beta=\sum_{k_1+\cdots k_n\ge
i-1}(\alpha_1^{k_1}\cdots\alpha_n^{k_n})\cdot
f_{k_1,\ldots,k_n}\in I_{i-1},\quad f_{k_1,\ldots,k_n}\in A,$$ by
using \eqref{2.6}, we see that $\beta\in I_i$ if and only if
\ga{2.7}
{\sum_{k_1+\cdots+k_n=i-1}(\alpha_1^{k_1}\cdots\alpha_j^{k_j-1}\cdots\alpha_n^{k_n})\cdot
k_jf_{k_1,\ldots,k_n}\,\,\in I_{i-1}} for all $1\le j\le n$. Since
$\{\alpha_1^{k_1}\cdots\alpha_n^{k_n}\,|\,k_1+\cdots+k_n\ge 1\}$
is a basis of $I_1$ locally and the lemma is true for $I_{i-1}$,
\eqref {2.7} is equivalent to \ga{2.8} {{\rm For} \ {\rm given} \
(k_1,\ldots,k_n) \ {\rm with} \ k_1+\cdots+k_n=i-1\\
k_jf_{k_1,\ldots,k_n}=0 \ \ {\rm for} \  {\rm all} \
j=1,\ldots,n\notag} which implies $f_{k_1,\ldots,k_n}=0$ whenever
$k_1+\cdots+k_n=i-1$. Thus $I_i$ is generated by
$\{\alpha_1^{k_1}\cdots\alpha_n^{k_n}\,\,|\, k_1+\cdots+k_n\ge
i\,\}$.
\end{proof}

\begin{lem}\label{lem2.3}
\begin{itemize}\item[(i)]$I_i=0$ when $i>n(p-1)$, and $\nabla(I_{i+1})\subset I_i\otimes\Omega^1_X$
for $i\ge 1$.
\item[(ii)]
$I_i/I_{i+1}\xrightarrow{\nabla}(I_{i-1}/I_i)\otimes\Omega^1_X$
are injective in the category of vector bundles for $1\le i\le
n(p-1)$. In particular, their composition \ga{2.9} {\nabla^i:
I_i/I_{i+1}\to (I_0/I_1)\otimes_{\sO_X}(\Omega^1_X)^{\otimes
i}=(\Omega^1_X)^{\otimes i}} is injective in the category of
vector bundles.
\end{itemize}
\end{lem}

\begin{proof} (i) follows from Lemma \ref{2.2} and Definition
\ref{2.1}. (ii) follows from \eqref {2.6}.

\end{proof}

In order to describe the image of $\nabla^i$ in \eqref{2.9}, we
recall a ${\rm GL}(n)$-representation $V^{[\ell]}\subset
V^{\otimes\ell}$ where $V$ is the standard representation of ${\rm
GL}(n)$. Let ${\rm S}_{\ell}$ be the symmetric group of $\ell$
elements with the action on $V^{\otimes\ell}$ by
$(v_1\otimes\cdots\otimes
v_{\ell})\cdot\sigma=v_{\sigma(1)}\otimes\cdots\otimes
v_{\sigma(n)}$ for $v_i\in V$ and $\sigma\in{\rm S}_{\ell}$. Let
$e_1,\,\ldots,\,e_n$ be a basis of $V$, for $k_i\ge 0$ with
$k_1+\cdots+k_n=\ell$ define \ga{2.10}
{v(k_1,\ldots,k_n)=\sum_{\sigma\in{\rm S}_{\ell}}(e_1^{\otimes
k_1}\otimes\cdots\otimes e_n^{\otimes k_n})\cdot\sigma }

\begin{defn}\label{defn2.4} Let $V^{[\ell]}\subset V^{\otimes\ell}$ be
the linear subspace generated by all vectors $v(k_1,\ldots,k_n)$
for all $k_i\ge 0$ satisfying $k_1+\cdots+k_n=\ell$. It is clearly
a sub-representation of ${\rm GL}(V)$. If $\sV$ is a vector bundle
of rank $n$,  the subbundle $\sV^{[\ell]}\subset
\sV^{\otimes\ell}$ is defined to be the associated bundle of the
frame bundle of $\sV$ (which is a principal ${\rm GL}(n)$-bundle)
through the representation $V^{[\ell]}$.
\end{defn}

In characteristic zero, $V^{[\ell]}$ is nothing but ${\rm
Sym}^{\ell}(V)$. When $char(k)=p>0$, we have $v(k_1,\ldots,k_n)=0$
if one of $k_1,\,\ldots,\,k_n$ is bigger than $p-1$. Thus
$V^{[\ell]}$ is in fact spanned by \ga{2.11}
{\{v(k_1,\ldots,k_n)\,|\, 0\le k_i\le p-1\,,\,
k_1+\cdots+k_n=\ell\,\}.} In general, $V^{[\ell]}$ is not
isomorphic to ${\rm Sym}^{\ell}(V)$, but it is easy to see
\ga{2.12} {V^{[\ell]}\cong{\rm Sym}^{\ell}(V) \quad {\rm when}
\quad 0<\ell<p\,.}

\begin{lem}\label{lem2.5} With the notation in Definition
\ref{defn2.4}, the composition \ga{2.13}
{\nabla^{\ell}:I_{\ell}/I_{\ell+1}\to (\Omega_X^1)^{\otimes\ell}}
of the $\sO_X$-morphisms in Lemma \ref{lem2.3} (ii) has image
$(\Omega^1_X)^{[\ell]}\subset (\Omega_X^1)^{\otimes\ell}$.
\end{lem}
\begin{proof} It is enough to prove the lemma locally. By Lemma
\ref{lem2.2}, $I_{\ell}/I_{\ell+1}$ is locally generated by
\ga{2.14} {\{\alpha_1^{k_1}\cdots\alpha_n^{k_n}\,\,|\,
k_1+\cdots+k_n= \ell\,\}.} By using formula \eqref{2.6}, we have
\ga{2.15}
{\nabla^{\ell}(\alpha_1^{k_1}\cdots\alpha_n^{k_n})=(-1)^{\ell}\sum_{\sigma\in{\rm
S}_{\ell}}({\rm d}x_1^{\otimes k_1}\otimes\cdots{\rm
d}x_n^{\otimes k_n})\cdot\sigma } which implies that
$\nabla^{\ell}(I_{\ell}/I_{\ell+1})=(\Omega^1_X)^{[\ell]}\subset
(\Omega_X^1)^{\otimes\ell}$.
\end{proof}

\begin{thm}\label{thm2.6} The filtration defined in Definition
\ref{defn2.1} is \ga{2.16} {0=V_{n(p-1)+1}\subset
V_{n(p-1)}\subset\cdots\subset V_1\subset V_0=V=F^*(F_*W)} which
has the following properties
\begin{itemize}\item[(i)]$\nabla(V_{i+1})\subset V_i\otimes\Omega^1_X$
for $i\ge 1$, and $V_0/V_1\cong W$.
\item[(ii)]
$V_i/V_{i+1}\xrightarrow{\nabla}(V_{i-1}/V_i)\otimes\Omega^1_X$
are injective morphisms of vector bundles for $1\le i\le n(p-1)$,
which induced isomorphisms $$\nabla^i: V_i/V_{i+1}\cong
W\otimes_{\sO_X}(\Omega^1_X)^{[i]},\quad 0\le i\le n(p-1).$$ In
particular, $V_i/V_{i+1}\cong W\otimes_{\sO_X}{\rm
Sym}^i(\Omega^1_X)$ for $i<p$.
\end{itemize}
\end{thm}

\begin{proof} It is a local problem to prove the theorem. Thus
$V_{n(p-1)+1}=0$ follows from Lemma \ref{lem2.2}, and (ii) follows
from Lemma \ref{lem2.3} and Lemma \ref{lem2.5}. (i) is nothing but
the definition.
\end{proof}

\bibliographystyle{plain}

\renewcommand\refname{References}

\end{document}